\documentclass[reqno]{amsart}

\usepackage{amsmath}
\usepackage{amssymb}
\usepackage{enumerate}

\newtheorem{thm}{Theorem}

\newtheorem{propn}{Proposition}

\newtheorem{nonumthm}{Theorem}

\theoremstyle{definition}

\newtheorem{eg}{Example}

\newcommand{\thmref}[1]{Theorem~\ref{#1}}
\newcommand{\propref}[1]{Proposition~\ref{#1}}

\newcommand{\secref}[1]{Section~\ref{#1}}
\newcommand{\agemo}{\mho}

\numberwithin{equation}{section}

\begin{document}
\baselineskip=18pt
\title {A Lie Algebra correspondence for a family of Finite
$p$-Groups}

\author{Paul J.\ Sanders}
\address
{Mathematics Institute\\ University of Warwick\\ Coventry\\ England CV4 7AL}
\email {pjs@maths.warwick.ac.uk}


\maketitle



\section{Introduction}

A conjecture of J.\ Moody\footnote{Private communication} asserts that for
each natural number $n$, a finite number of primes can be excluded so that if
$p$ is not one of the excluded primes, there is a (non-natural) 
one-to-one correspondence between isomorphism classes of groups of order $p^n$ 
and isomorphism classes of nilpotent $p^n$-element
${\mathbb F}_p[T]/(T^n)$-Lie algebras.
In this paper we consider a special case of this conjecture. Specifically, we 
show

\begin{thm}
\label{maintheorem}
Let $n$ be a natural number and $p$ a prime with $p\geq n$. Then there exists a
one-to-one correspondence between isomorphism classes of groups
$P$ of order $p^n$ in which $\agemo_1([P,P]) = 1$, and isomorphism classes of
nilpotent $p^n$-element ${\mathbb F}_p[T]/(T^n)$-Lie algebras $L$ in which
$T[L,L] = 0$.
\end{thm}

\noindent The condition on $p$ enables us to regard the groups in the 
statement as certain Lie rings of order $p^n$ (by means of the 
Magnus-Lazard correspondence, the main properties of which are summarised in
\secref{proofs}). We prove this 
theorem by establishing a 
correspondence between this category and the category of Lie algebras in the
statement.

When $p\geq n$, any $p^n$-element nilpotent ${\mathbb F}_p[T]/(T^n)$-Lie 
algebra $L$ has an associated group of order $p^n$ and exponent $p$ 
(from the applying the Magnus-Lazard correspondence to the underlying 
${\mathbb F}_p$-Lie algebra). Using this observation, \thmref{maintheorem}
can be used to give an expression for the number of groups of order $p^n$ whose
derived subgroup has exponent dividing $p$. Denoting by ${\mathcal E}_p^n$,
a transversal for the isomorphism classes of groups of order $p^n$ and
exponent $p$, and for such a group $E$, letting $N_E$ denote the endomorphisms
$\sigma : E \rightarrow E$ with the property that 
${\rm Im}\, \sigma \subseteq Z(E)$ and $\sigma^n = e_E$ (the endomorphism of
$E$ sending the entire group to the identity), we will obtain 

\begin{thm}
\label{expression}
Let $n$ be a natural number and $p$ a prime with $p\geq n$. Then the number of
isomorphism classes of groups of order $p^n$ whose derived subgroup has
exponent dividing $p$, is given by
\begin{equation*}
\sum_{E\in{\mathcal E}_p^n} \sum_{\phi\in N_E} 
\frac{1}{\vert\phi^{{\rm Aut}(E)}\vert},
\end{equation*}
where ${\rm Aut}(E)$ acts on ${\mathcal E}_p^n$ by conjugation.
\end{thm}

This paper comprises part of the author's doctoral dissertation 
\cite{mythesis}. The author would like to thank his advisor J.\ Moody for
outlining the conjecture referred to above.

\section{A Correspondence between certain non-associative Rings and Algebras}
If $C$ is any commutative ring with 1, we define a $C$-algebra to be a unital 
$C$-module $A$ equipped with an element of ${\rm Hom}_C(A\otimes_C A,A)$
giving the multiplication in $A$ (multiplication will be denoted by a bracket
$[.,.]$). In this context, a ring will be taken to mean a 
${\mathbb Z}$-algebra.

In this section we establish, for an arbitrary prime $p$ and natural number 
$n$, a one-to-one correspondence between isomorphism classes in the category
${\mathcal R}_p^n$ of rings $R$ of order $p^n$ in which $p[R,R] = 0$, and
isomorphism classes in the category ${\mathcal A}_p^n$ of $p^n$-element
${\mathbb F}_p[T]/(T^n)$-algebras $A$ in which $T[A,A] = 0$. We use the type
invariants of the underlying modules to divide up the correspondence.

\begin{thm}
\label{gencorr}
Let $p$ be a prime, n be a natural number, and suppose that 
$\lambda_1 \geq \cdots \geq \lambda_t$ is a partition of $n$. Let $U$ be
the Abelian $p$-group 
${\mathbb Z}/p^{\lambda_1}{\mathbb Z}\oplus\cdots\oplus
{\mathbb Z}/p^{\lambda_t}{\mathbb Z}$, and $V$ be the 
${\mathbb F}_p[T]/(T^n)$-module
${\mathbb F}_p[T]/(T^{\lambda_1})\oplus\cdots\oplus
 {\mathbb F}_p[T]/(T^{\lambda_t})$.
Then 
\begin{enumerate}
\item There exists a one-to-one correspondence between isomorphism classes of
rings defined on $U$ which belong to ${\mathcal R}_p^n$, and isomorphism
classes of ${\mathbb F}_p[T]/(T^n)$-algebras defined on $V$ which belong to
${\mathcal A}_p^n$.
\item This correspondence sends isomorphism classes of Lie rings onto
isomorphism classes of ${\mathbb F}_p[T]/(T^n)$-Lie algebras, and for such
classes, the property of nilpotency {\rm(}when it exists{\rm)} is preserved.
\end{enumerate}
\end{thm}

\begin{proof}
Let $\Omega_1(U)=\{x\in U: px=0\}, \agemo_1(U)=\{px : x\in U\}, 
\Omega_1(V)=\{x\in V:Tx=0\}, \agemo_1(V)=\{Tx:x\in V\}$. Denote the quotient
${\mathbb Z}$-module $U/\agemo_1(U)$ by $\overline U$, the quotient
${\mathbb F}_p[T]$-module $V/\agemo_1(V)$ by $\overline V$, and observe that
$T$-action on $\Omega_1(V)$ and $\overline V$ is trivial. 
Then any ring 
structure $R$ on $U$ with $p[R,R]=0$ is uniquely determined by an element of 
${\rm Hom}_{\mathbb Z}(\overline U \otimes_{\mathbb Z} 
\overline U,\Omega_1(U))$, and any ${\mathbb F}_p[T]$-algebra structure $A$
on $V$ with $T[A,A] = 0$ is uniquely determined by an element of
${\rm Hom}_{{\mathbb F}_p[T]}(\overline V\otimes_{{\mathbb F}_p[T]}
\overline V,\Omega_1(V))$. Moreover, the isomorphism classes of ring
structures $R$ on $U$ with $p[R,R]=0$ are in one-to-one correspondence with
the induced action of ${\rm Aut}_{\mathbb Z}(U)$ on 
${\rm Hom}_{\mathbb Z}(\overline U \otimes_{\mathbb Z}\overline U,
\Omega_1(U))$, and the isomorphism classes of ${\mathbb F}_p[T]$-algebra
structures $A$ on $V$ with $T[A,A]=0$ are in one-to-one correspondence with
the induced action of ${\rm Aut}_{{\mathbb F}_p[T]}(V)$ on
${\rm Hom}_{{\mathbb F}_p[T]}(\overline V\otimes_{{\mathbb F}_p[T]}
\overline V,\Omega_1(V))$.
 
Now fix a ${\mathbb Z}$-basis $(u_1,\ldots,u_t)$
of $U$ and an ${\mathbb F}_p[T]$-basis $(v_1,\ldots,v_t)$ of $V$, both
corresponding to the type invariants $(\lambda_1,\ldots,\lambda_t)$. Then we
have ${\mathbb F}_p$-bases 
$\{\overline{u_i}\otimes\overline{u_j} : 1\leq i,j\leq t\}$,
$\{\overline{u_i}\otimes\overline{u_j} : 1\leq i,j\leq t\}$,
$\{p^{\lambda_1-1}u_1,\ldots,p^{\lambda_t-1}u_t\}$, 
$\{T^{\lambda_1-1}v_1,\ldots,T^{\lambda_t-1}v_t\}$ of
$\overline U \otimes_{\mathbb Z}\overline U$,
$\overline V\otimes_{{\mathbb F}_p[T]}\overline V$,
$\Omega_1(U)$ and $\Omega_1(V)$ respectively. These give 
${\mathbb F}_p$-vector space isomorphisms $\psi,\nu$
\begin{equation*}
\begin{array}{rrccc}
\psi&:&\overline U \otimes_{\mathbb Z}\overline U&\longrightarrow&
\overline V\otimes_{{\mathbb F}_p[T]}\overline V\\
&&\overline{u_i}\otimes\overline{u_j}&\longmapsto&\overline{v_i}\otimes
\overline{v_j}
\end{array}
\qquad
\begin{array}{rrccc}
\nu&:&\Omega_1(U)&\longrightarrow&\Omega_1(V)\\
&&p^{\lambda_i-1}u_i&\longmapsto&T^{\lambda_i-1}v_i
\end{array}
\end{equation*}
and hence an ${\mathbb F}_p$-vector space isomorphism
\begin{equation*}
\theta : {\rm Hom}_{\mathbb Z}(\overline U \otimes_{\mathbb Z}\overline U,
\Omega_1(U)) \longrightarrow {\rm Hom}_{{\mathbb F}_p[T]}(\overline V
\otimes_{{\mathbb F}_p[T]}\overline V,\Omega_1(V))
\end{equation*}
With respect to the above bases, we have linear actions of
${\rm GL}_t({\mathbb F}_p) \times {\rm GL}_t({\mathbb F}_p)$ on
${\rm Hom}_{\mathbb Z}(\overline U \otimes_{\mathbb Z}\overline U,\Omega_1(U))$
and
${\rm Hom}_{{\mathbb F}_p[T]}(\overline V\otimes_{{\mathbb F}_p[T]}\overline V,
\Omega_1(V))$, and $\theta$ is equivariant for these actions. Moreover, the
induced actions of ${\rm Aut}_{\mathbb Z}(U)$ and 
${\rm Aut}_{{\mathbb F}_p[T]}(V)$ on 
${\rm Hom}_{\mathbb Z}(\overline U \otimes_{\mathbb Z}\overline U,\Omega_1(U))$
and 
${\rm Hom}_{{\mathbb F}_p[T]}(\overline V\otimes_{{\mathbb F}_p[T]}\overline V,
\Omega_1(V))$ respectively, both factor through the action of
${\rm GL}_t({\mathbb F}_p) \times {\rm GL}_t({\mathbb F}_p)$. To prove part 1 
of the theorem therefore, it suffices to show that they factor through the same
subgroup of ${\rm GL}_t({\mathbb F}_p) \times {\rm GL}_t({\mathbb F}_p)$.

First, define integers $1\leq i_1 < i_2 < \cdots < i_r \leq t$ such that
$\{\lambda_{i_1},\lambda_{i_2},\ldots,\lambda_{i_r}\} =
\{\lambda_1,\ldots,\lambda_t\}$ and $\lambda_{i_j} > \lambda_{i_{j+1}}$
for $1\leq j < r$. Now for each $j=1,\ldots,r$, define $d_j$ to be the number
of $\lambda_l$'s equal to $\lambda_{i_j}$. Then if 
$\pi\in{\rm Aut}_{\mathbb Z}(U)$, we see that relative to the basis
$\{p^{\lambda_1-1}u_1,\ldots,p^{\lambda_t-1}u_t\}$ of $\Omega_1(U)$,
$\pi\vert_{\Omega_1(U)}$ is represented by a matrix of the form
\begin{equation}
\label{restrict}
\begin{pmatrix}
M_1&0&\cdots&0&0\\
*&M_2&\cdots&0&0\\
\vdots&\vdots&\ddots&\vdots&\vdots\\
*&*&\cdots&M_{r-1}&0\\
*&*&\cdots&*&M_r\\
\end{pmatrix}
\in{\rm GL}_t({\mathbb F}p)
\end{equation}
where $M_j\in{\rm GL}_{d_j}({\mathbb F}_p)$ for each $j=1,\ldots,r$.
Moreover, relative to the basis $(\overline{u_1},\ldots,\overline{u_t})$ of
$\overline U$, the element of ${\rm Aut}_{\mathbb Z}(\overline U)$ is
represented by a matrix of the form
\begin{equation}
\label{induce}
\begin{pmatrix}
M_1&*&\cdots&*&*\\
0&M_2&\cdots&*&*\\
\vdots&\vdots&\ddots&\vdots&\vdots\\
0&0&\cdots&M_{r-1}&*\\
0&0&\cdots&0&M_r\\
\end{pmatrix}
\in{\rm GL}_t({\mathbb F}_p)
\end{equation}
where the $M_j$ are the same as the diagonal blocks appearing in the
representation \eqref{restrict} for $\pi\vert_{\Omega_1(U)}$. Conversely,
given matrices $E,F\in{\rm GL}_t({\mathbb F}_p[T])$ where $E$ is of the form
\eqref{restrict}, $F$ is of the form \eqref{induce} and the diagonal blocks
of $E$ and $F$ coincide, we can lift $F$ to an element of 
${\rm Aut}_{\mathbb Z}(U)$ whose restriction to $\Omega_1(U)$ relative to
$\{p^{\lambda_1-1}u_1,\ldots,p^{\lambda_t-1}u_t\}$ is $E$. 

It is straightforward to see that a similar situation holds for
${\rm Aut}_{{\mathbb F}_p[T]}(V)$, i.e.\ that the subgroup of 
${\rm GL}_t({\mathbb F}_p)\times{\rm GL}_t({\mathbb F}_p)$ consists of
matrices $(E,F)$ where $E$ is of the form \eqref{induce}, $F$ is of the
form \eqref{restrict}, and the diagonal blocks of $E$ and $F$ coincide.
Therefore $\theta$ induces a one-to-one correspondence between isomorphism 
classes. 

Maintaining the above notation, we now show part 2 of the theorem. So let
$\epsilon_R\in{\rm Hom}_{\mathbb Z}(\overline U \otimes_{\mathbb Z} 
\overline U,\Omega_1(U))$ and 
$\gamma_A\in{\rm Hom}_{{\mathbb F}_p[T]}(\overline V
\otimes_{{\mathbb F}_p[T]} \overline V, \Omega_1(V))$ with $\epsilon_R\theta =
\gamma_A$. For any $1\leq i,j \leq t$, it follows by the construction of
$\theta$ that 
$(\overline{u_i}\otimes\overline{u_j})\epsilon_R =
-(\overline{u_j}\otimes\overline{u_i})\epsilon_R$ if and only if
$(\overline{v_i}\otimes\overline{v_j})\gamma_A =
-(\overline{v_j}\otimes\overline{v_i})\gamma_A$, and
$(\overline{u_i}\otimes\overline{u_j})\epsilon_R = 0$ if and only if
$(\overline{v_j}\otimes\overline{v_i})\gamma_A = 0$. Therefore,
$\epsilon_R$ factors through $\bigwedge_2(\overline U)$ if and only if 
$\gamma_A$ factors through $\bigwedge_2(\overline V)$. 

We now show that the Jacobi identity
holds in $R$ if and only if it holds in $A$, so for elements $x,y,z\in R$
and $a,b,c\in A$, let $J_R(x,y,z)=[[x,y],z]+[[y,z],x]+[[z,x],y]$ and
$J_A(a,b,c)=[[a,b],c]+[[b,c],a]+[[c,a],b]$. Observe that the Jacobi identity
holds in $R$ if and only if $J_R(u_i,u_j,u_k)=0$ for all 
$i,j,k=1,\ldots,t$, and the Jacobi identity holds in $A$ if and only if
$J_A(v_i,v_j,v_k)=0$ for all $i,j,k=1,\ldots,t$. Now if $\lambda_t\geq2$,
we have $\Omega_1(U)\subseteq\agemo_1(U)$, and so for any $x,y,z\in R$ we have
$[[x,y],z] = 0$. Similarly, the Jacobi identity holds in $A$. Therefore we
may assume that $\lambda_t=1$ and so $d_r$ (defined above) equals the number of
type invariants equal to 1. For any $1\leq i,j\leq t$ choose integers 
$\{\alpha_{ij}^l\}_{l=1}^t$ such that
\begin{equation}
\label{structure}
[u_i,u_j]=\sum_{l=1}^t \alpha_{ij}^lp^{\lambda_l-1}u_l\qquad{\rm and}\qquad
[v_i,v_j]=\sum_{l=1}^t \alpha_{ij}^lT^{\lambda_l-1}v_l.
\end{equation}
(such integers exist since we are assuming $\epsilon_R\theta=\gamma_A$).
Now since $\lambda_e - 1 \geq 1$ for $1\leq e\leq t-d_r$, and $\lambda_e = 1$
for $t-d_r+1\leq e\leq t$, it follows that for $1\leq i,j,k \leq t$ we have 
\begin{eqnarray}
\label{JacobiR}
[[u_i,u_j],u_k] = \sum_{l=t-d_r+1}^t\alpha_{ij}^l[u_l,u_k]&=&
\sum_{l=t-d_r+1}^t\alpha_{ij}^l\left(\sum_{m=1}^t\alpha_{lk}^mp^{\lambda_m-1}
u_m\right)\nonumber\\
&=&\sum_{m=1}^t\left(\sum_{l=t-d_r+1}^t\alpha_{ij}^l\alpha_{lk}^m\right)
p^{\lambda_m-1}u_m
\end{eqnarray}
and
\begin{equation}
\label{JacobiA}
[[v_i,v_j],v_k] = \sum_{m=1}^t\left(\sum_{l=t-d_r+1}^t\alpha_{ij}^l
\alpha_{lk}^m\right)T^{\lambda_m-1}v_m.
\end{equation}
By permuting $i,j,k$ cyclically in \eqref{JacobiR} and \eqref{JacobiA}, we
obtain
\begin{equation*}
J_R(u_i,u_j,u_k)=\sum_{m=1}^t\left(\sum_{l=t-d_r+1}^t\alpha_{ij}^l
\alpha_{lk}^m + \alpha_{ki}^l\alpha_{lj}^m + \alpha_{jk}^l\alpha_{li}^m
\right)p^{\lambda_m-1}u_m
\end{equation*}
and
\begin{equation*}
J_A(v_i,v_j,v_k)=\sum_{m=1}^t\left(\sum_{l=t-d_r+1}^t\alpha_{ij}^l
\alpha_{lk}^m + \alpha_{ki}^l\alpha_{lj}^m + \alpha_{jk}^l\alpha_{li}^m
\right)T^{\lambda_m-1}v_m,
\end{equation*}
from which it follows that $J_R(u_i,u_j,u_k)=0$ if and only if
$J_A(v_i,v_j,v_k)=0$. Hence $\epsilon_R$ defines a Lie ring structure on $U$
if and only if $\gamma_A$ defines an ${\mathbb F}_p[T]/(T^n)$-Lie algebra
structure on $A$. 

To complete the proof of part 2, let 
$\epsilon_R\in{\rm Hom}_{\mathbb Z}(\overline U \otimes_{\mathbb Z} 
\overline U,\Omega_1(U))$ where $R$ is a Lie ring on $U$, and write 
$\gamma_A=\epsilon_R\theta$ so that $A$ is an ${\mathbb F}_p[T]/(T^n)$-Lie 
algebra on $V$. Then $R/\agemo_1(R) \cong A/\agemo_1(A)$ as ${\mathbb F}_p$-Lie
algebras (to see this, consider the structure constants given by 
\eqref{structure}). Now since $[\agemo_1(R),R]=0$ it follows that nilpotency of
$R$ is equivalent to nilpotency of $R/\agemo_1(R)$. Similarly, nilpotency of
$A$ is equivalent to nilpotency of $A/\agemo_1(A)$, Hence $R$ is nilpotent if
and only if $A$ is nilpotent. 
\end{proof}

\section {Proofs of the Main Theorems}
\label{proofs}
As mentioned in the introduction, the proofs of \thmref{maintheorem} and
\thmref{expression} utilise a Lie ring correspondence first discovered
by Magnus in \cite{magnus}, and then later, independently, by Lazard in
\cite{lazard}. For a prime
$p$, let $\Gamma_p$ denote the category of finite $p$-groups whose nilpotency
class is less than $p$, and let $\Lambda_p$ denote the category of finite
nilpotent Lie rings whose order is a power of $p$ and whose nilpotency class
is less than $p$. The main properties of the correspondence are summarised in
the following theorem.
 
\begin{nonumthm}[Magnus \cite{magnus}, Lazard \cite{lazard}]
Let $p$ be a prime, $P$ a group in $\Gamma_p$ and $L$ a Lie ring in $\Lambda_p$.
Then there exists a Lie ring ${\mathcal L}_p(P)$ in $\Lambda_p$ and a $p$-group
${\mathcal G}_p(L)$ in $\Gamma_p$ such that
\begin{enumerate}
\item $L$ has the same underlying set as ${\mathcal G}_p(L)$.
\item $P$ has the same underlying set as ${\mathcal L}_p(P)$.
\item The identity of $P$ is the identity of the underlying Abelian group of
${\mathcal L}_p(P)$. 
\item\label{inversea} The Lie ring ${\mathcal L}_p({\mathcal G}_p(L))$
coincides with $L$.
\item\label{inverseb} The $p$-group ${\mathcal G}_p({\mathcal L}_p(P))$
coincides with $P$.
\item\label{covariant} ${\mathcal G}_p$ and ${\mathcal L}_p$ are covariant
functors between the respective categories (the morphisms being the obvious
structure preserving maps). 
\item\label{orderscoincide} The order of an element $g\in P$ coincides with its
additive order in
${\mathcal L}_p(P)$.
\item For $x,y\in P$, the group commutator $[x,y]_P$ equals 1 if and only if
the Lie bracket $[x,y]_{{\mathcal L}_p(P)}$ equals 0.
\item Subgroups of $P$ coincide with Lie-subrings of ${\mathcal L}_p(P)$, and
normal subgroups of $P$ coincide with ideals of ${\mathcal L}_p(P)$.
\item\label{commcoincide} For normal subgroups $H$ and $K$ of $P$, the 
commutator group $[H,K]_P$
coincides with the ideal of ${\mathcal L}_p(P)$ generated by all Lie brackets
of elements of $H$ with elements of $K$.
\end{enumerate}
\end{nonumthm}

The functors ${\mathcal G}_p$ and ${\mathcal L}_p$ arise from the 
Baker-Campbell-Hausdorff formula and its ``inversion''. The reader is referred 
to \cite{lazard} or \cite{mythesis} for the exact details --- we need only the 
properties asserted in the above theorem.

By Hall's commutator collecting process in \cite{hall}, any $p$-group in 
$\Gamma_p$ is regular. In such a $p$-group $P$, the subset 
$\Omega_i(P)=\{x\in P: x^{p^i}=1\}$ is a subgroup, and denoting the exponent of
$P$ by $p^{\mu(P)}$, we obtain a partition 
$(\omega_1(P),\ldots,\omega_{\mu(P)}(P))$ of $n$ by defining
$p^{\omega_i(P)} = \vert \Omega_i(P) : \Omega_{i-1}(P) \vert$. For
$1\leq i\leq\omega_1(P)$, we define $\mu_i(P)$ to be the number of $\omega_j$'s
greater than or equal to $i$, and this gives us a dual partition of $n$. We
say that $P$ is of type $(\mu_1(P),\ldots,\mu_{\omega_1(P)}(P))$. Observe that
since the above Lie ring correspondence preserves the order of an element it
follows that for a $p$-group $P$ in $\Gamma_p$, the underlying Abelian 
$p$-group of ${\mathcal L}_p(P)$ is also of type 
$(\mu_1(P),\ldots,\mu_{\omega_1(P)}(P)$.

\begin{propn}
\label{regcorr}
Let $n$ be a natural number and $p$ a prime with $p\geq n$. Then for any
partition $\lambda_1\geq\cdots\geq\lambda_t$ of $n$, there is a one-to-one
correspondence between isomorphism classes of regular $p$-groups $P$ of 
type $(\lambda_1,\ldots,\lambda_t)$ in which $\agemo_1([P,P])=1$, and 
isomorphism classes of nilpotent ${\mathbb F}_p[T]/(T^n)$-
Lie algebras $L$ of type $(\lambda_1,\ldots,\lambda_t)$ in which $T[L,L]=0$.
\end{propn}

\begin{proof}
First observe that since $p\geq n$, any group $P$ of order $p^n$ belongs to 
$\Gamma_p$ and any nilpotent Lie ring of order $p^n$ belongs to $\Lambda_p$.
Moreover, by property \ref{commcoincide} of the Magnus-Lazard correspondence, 
the derived subgroup of $P$ coincides with the derived subring of 
${\mathcal L}_p(P)$, and so by property \ref{orderscoincide}, the exponent of 
$[P,P]$ coincides with the additive exponent of
$[{\mathcal L}_p(P),{\mathcal L}_p(P)]$. Therefore there is a one-to-one
correspondence between the $p$-groups in the statement and isomorphism classes
of nilpotent $p^n$-element Lie rings $R$ of type $(\lambda_1,\ldots,\lambda_t)$
in which $p[R,R]=0$. The result now follows from \thmref{gencorr}.
\end{proof} 

\noindent We now have
\begin{proof}[Proof of \thmref{maintheorem}]
Apply \propref{regcorr} to each partition of $n$. 
\end{proof}

By ignoring the $T$-action on a nilpotent $p^n$-element 
${\mathbb F}_p[T]/(T^n)$-Lie algebra, we obtain a nilpotent $n$-dimensional
${\mathbb F}_p$-Lie algebra $K$. Moreover, the isomophism classes of nilpotent 
${\mathbb F}_p[T]/(T^n)$-Lie algebra structures $L$ on $K$ in which $T[L,L]=0$,
correspond to the
orbits of the conjugation action of ${\rm Aut}_{\rm Lie}(K)$ on the set $C(K)$
defined as
\begin{equation*}
\{\sigma\in{\rm Hom}_{{\mathbb F}_p}(K,K) : \sigma^n=0,\
{\rm Im}\,\sigma\subseteq Z(K),\
[x\sigma,y] = [x,y]\sigma\ {\rm for\ all}\ x,y\in K\}
\end{equation*}

For any $\sigma\in C(K)$, the type invariants of the corresponding
${\mathbb F}_p[T]/(T^n)$-module are the dimensions of the Jordan blocks of
$\sigma$. So for a partition $\underline \lambda$ of $n$, if we denote by
$C(K,\underline \lambda)$ the elements $\sigma\in C(K)$ with Jordan block 
dimensions 
$\underline \lambda$, then the number of isomorphism classes of
nilpotent ${\mathbb F}_p[T]/(T^n)$-Lie algebra structures $L$ on $K$ with
$T[L,L]=0$ equals 
$\sum_{\sigma\in C(K,\underline\lambda)}1/\vert \sigma^{{\rm Aut}_{\rm Lie}(K)}
\vert$. 

\begin{propn}
\label{liesum}
Let $p\geq n$ and denote by ${\mathcal K}_p^n$ a transversal for the isomorphism
classes of nilpotent $n$-dimensional ${\mathbb F}_p$-Lie algebras. Then 
for a partition 
$\underline \lambda$ of $n$, the
number of isomorphism classes of $p$-groups $P$ of type
$\underline \lambda$ in which $\agemo_1([P,P])=1$ is given by
\begin{equation*}
\sum_{K\in{\mathcal K}_p^n}\sum_{\sigma\in C(K,\underline \lambda)}
\frac{1}{\vert \sigma^{{\rm Aut}_{\rm Lie}(K)} \vert}.
\end{equation*}
\end{propn}
\begin{proof}
Follows immediately from the preceding discussion and \propref{regcorr}.
\end{proof}

\begin{eg}
For $p\geq 7$, the 7-dimensional nilpotent ${\mathbb F}_p$-Lie algebras are
known from applying the Magnus-Lazard correspondence to Wilkinson's
classification \cite{wilkinson} of group of order $p^7$ and exponent $p$. This
proposition gives an expression for the $p$-groups of types $(2,1,1,1,1,1)$
and $(2,2,1,1,1)$ in terms of these Lie algebras.
\end{eg}

\begin{proof}[Proof of \thmref{expression}]
Fixing $p\geq n$, the Magnus-Lazard theorem gives a one-to-one correspondence
between isomorphism classes of groups of order $p^n$ with exponent $p$, and 
isomorphism classes of nilpotent $p^n$-element Lie rings with 
additive exponent $p$. Therefore, we can take ${\mathcal K}_p^n =
\{{\mathcal L}_p(E) : E\in{\mathcal E}_p^n\}$. Now if $E\in{\mathcal E}_p^n$ 
and $\sigma\in C({\mathcal L}_p(E))$ we see easily from the Magnus-Lazard 
theorem that $\sigma\in N_E$. Conversely, any element $\tau\in N_E$ is a 
nilpotent endomorphism of ${\mathcal L}_p(E)$ with ${\rm Im}\,\tau\subseteq
Z({\mathcal L}_p(E))$, and for any $x,y\in{\mathcal L}_p(E)$ we have
$[x,y]\tau=[x\tau,y\tau]=0=[x\tau,y]$. Therefore, $N_E = C({\mathcal L}_p(E)$, 
and so the result follows immediately after summing the expression in 
\propref{liesum} over all partitions of $n$.
\end{proof}

\bibliography{connect}
\bibliographystyle{plain}
 
\end{document}